\documentclass{amsart}
\usepackage{amsmath, amssymb}
\usepackage
{hyperref}
\hypersetup{colorlinks=true,citecolor=blue,linkcolor=blue,urlcolor=blue,
pdfstartview=FitH }
\begin{document}

\title{The relative trace formula in analytic number theory}
\date{}

\author{Valentin Blomer}
\address{Mathematisches Institut, Endenicher Allee 60, 53115, Bonn, Germany}\email{blomer@math.uni-bonn.de} 

\begin{abstract}  We discuss a variety of applications of the relative trace formula of Kuznetsov type in analytic number theory and the theory of automorphic forms.  \end{abstract}



\thanks{Author partially supported by DFG grant BL 915/2-2}

\maketitle

\section{The Poisson summation formula}

The mother of all trace formulae is the Poisson summation formula. If $f : \Bbb{R} \rightarrow \Bbb{C}$ is a sufficiently nice function, say Schwartz class for simplicity, then 
\begin{equation}\label{poisson}
\sum_{n \in \Bbb{Z}} f(n)= \sum_{n \in \Bbb{Z}} \widehat{f}(n)
\end{equation}
where
$$\widehat{f}(y) = \int_{-\infty}^{\infty} f(x) e(-xy) \, {\rm d}x, \quad e(x) = e^{2\pi ix}$$
is the Fourier transform. Applying \eqref{poisson} formally with $f(x) = x^{-s}$, one obtains the functional equation of the Riemann zeta function, and this heuristic argument can easily be made rigorous. More generally, the Poisson summation formula implies  the functional equations for all Dirichlet $L$-functions \cite[Section 4.6]{IK}, and is in fact essentially a re-statement of those. 

Why is the Poisson summation formula a trace formula? Given a function $f$ as above, we consider the   convolution operator
$$L_f : L^2(\Bbb{R}/\Bbb{Z}) \rightarrow L^2(\Bbb{R}/\Bbb{Z})$$
given by 
$$ L_f(g)(x)  =    \int_{\Bbb{R}} f(x-y) g(y) \, {\rm d}y  = \int_{\Bbb{R}/\Bbb{Z}}k(x, y)g(y)\, {\rm d}y$$
where
\begin{equation}\label{ker}
    k(x, y) = \sum_{n \in \Bbb{Z}} f(x-y+n).
    \end{equation}
 The functions $e_n(x) = e(nx)$ are an orthonormal basis of eigenfunctions: $L_f(e_n) = \widehat{f}(n) e_n$. Computing the trace of $L_f$ in two ways we obtain
 $$\sum_{n\in \Bbb{Z}} \widehat{f}(n) = \int_{\Bbb{R}/\Bbb{Z}} k(x, x) \, {\rm d}x = \sum_{n \in \Bbb{Z}} f(n).$$
An analogous formula can be established for the multi-dimensional torus $(\Bbb{R}/\Bbb{Z})^n$. The Poisson summation formula is one of the most frequently used tools in analytic number theory, and we give two standard applications:

a) \emph{The P\'olya-Vinogradov inequality} (\cite{Po, Vi}, see \cite[Section 12.4]{IK} for a modern treatment). If $\chi$ is a primitive Dirichlet  character modulo $q$, then
$$\sum_{n \leq x} \chi(n) \ll q^{1/2} \log q$$
for all $x \geq 1$.  
This shows cancellation in character sums of length slightly larger than the square-root of the modulus. The proof is simple:  we split the sum over $n$ into residue classes modulo $q$, smooth the characteristic function on $[1, x]$ a little bit and apply the Poisson summation formula. This leads to the bound
$$ \sum_{n \leq x} \chi(n)  \ll 1 + \frac{1}{\sqrt{q}}\Big| \sum_{n \in \Bbb{Z}} \bar{\chi}(n) w(n)\Big|$$
where $w(n) \ll \min(q/|n|, q^2/|n|^2)$. 
The zero frequency $n=0$ vanishes, and a trivial estimate on the non-zero frequencies gives the desired bound.

b) \emph{The number of lattice points in growing discs} (\cite{Sz}, see \cite[Section 4.4]{IK} or \cite[Section 10.2.6]{Cohen} for a modern treatment with complete details). Given a disc of radius $\sqrt{R}$ centered at the origin, it contains about $\pi R$ lattice points as $R \rightarrow \infty$. The error term can easily be bounded $O(R^{1/2})$ (the circumference of the disc), but an application of the two-dimensional Poisson summation formula improves this error to $O(R^{1/3})$. Here the integral transform involves a Bessel function $J_0$ as an integral kernel.  This can also be interpreted as an instance of Weyl's law on the flat torus. 

The experienced reader will observe that both applications are essentially nothing but applications of functional equations of suitable zeta functions.

\section{Harmonic analysis on upper half plane}\label{sec2}

The group $(\Bbb{R}/\Bbb{Z})^n$ is abelian, and for the rest of this note we will consider spectral summation formulae for non-abelian groups. 
 In this section the underlying group is $G = {\rm SL}_2(\Bbb{R})$ and for simplicity we fix the lattice $\Gamma = {\rm SL}_2(\Bbb{Z})$. For $K = {\rm SO}(2)$, the quotient is $X := \Gamma \backslash G/K$ is in natural bijection with the quotient of the upper half plane $\Bbb{H} = \{z \in \Bbb{C} \mid \Im z > 0\}$ by $\Gamma$. The space $L^2(X)$  is acted on by various self-adjoint and pairwise commuting operators: the hyperbolic Laplace operator
$$\Delta = -y^2(\partial_x^2 + \partial_y^2),$$
its non-archimedean counterparts, the Hecke operators
$$T_n(\phi)(z) = \frac{1}{\sqrt{n}} \sum_{ad = n} \sum_{b \, (\text{mod } d)}   \phi\left(\frac{az+b}{d}\right)  $$
for   $n\in \Bbb{N}$, and the reflection operator $T_{-1} (\phi)(z) = \phi(-\bar{z}).$ 
Since $X$ is not compact, we have a spectral decomposition
$$L^2(X) = L^2_{\text{disc}}(X) \oplus L^2_{\text{cont}}(X) =    L^2_{\text{cusp}}(X) \oplus L^2_{\text{res}}(X) \oplus L^2_{\text{cont}}(X)$$
where $L^2_{\text{cusp}}(X)$ is the   subspace of $\phi \in L^2(X)$ satisfying $\int_0^1 \phi(x + iy) \, {\rm d}x = 0$ for all $y > 0$,    $L^2_{\text{cont}}(X)$ is the continuous spectrum spanned by incomplete Eisenstein series, and in our case $L^2_{\text{disc}}(X)$ consists only of the one-dimensional space of constant functions.  A non-constant joint eigenfunction of $\Delta$,  all $T_n$ ($n \in \Bbb{N}$) and $T_{-1}$ with respective eigenvalues $\lambda_{\phi} = t_{\phi}^2 + 1/4$,  $\lambda_{\phi}(n)$ with $n\in \Bbb{N}$, and $\epsilon_{\phi} \in \{\pm 1\}$ has a Fourier expansion
$${\phi}(x+iy) = \text{constant term } +  \sum_{n \not = 0} \frac{\lambda_{\phi}(n)}{\sqrt{|n|}} e(nx) W_{it_{\phi}}(4\pi |n|y)$$
where $W_{it}(y) = W_{0, it}(y)$ is the standard Whittaker function and for $n < 0$ we define 
$\lambda_{\phi}(n) =  (-1)^{(1 - \epsilon_{\phi})/2} \lambda_{\phi}(-n).$ 
 By definition, the constant term vanishes if and only if ${\phi} \in L^2_{\text{cusp}}(X)$. The fact that the Hecke eigenvalues are essentially the Fourier coefficients of ${\phi}$ is a specialty of the group ${\rm GL}(2)$, this is not the case for general groups. 

We proceed to derive a sort of  analogue of the Poisson summation formula. Given a sufficiently nice function $f : K\backslash G/K \rightarrow \Bbb{C}$, say smooth and compactly supported, we construct a kernel function $K : X \times X \rightarrow \Bbb{C}$ by
\begin{equation}\label{defK}
K(x_1, x_2) = \sum_{\gamma \in \Gamma} f(x_1^{-1} \gamma x_2),
\end{equation}
which is is the analogue of \eqref{ker} and  gives rise to a convolution operator $L_f(g)(x_1) = \int_{X} K(x_1, x_2) g(x_2) \, {\rm d}x_2.$ Again we try to compute the trace of $L_f$ in two ways: as the sum over the eigenvalues and as the integral of the kernel over the diagonal. As our space $X$ is non-compact, there is a convergence problem that is usually solved by applying a suitable truncation operator. The integral of the kernel is computed explicitly by splitting the sum over $\gamma \in \Gamma$ into conjugacy classes. In this way one obtains (the most basic form of) the Selberg trace formula as a sum over the spectrum on the one side and a sum over conjugacy classes of matrices $\gamma \in \Gamma$ on the other, which features in particular class numbers and regulators of quadratic number fields in its hyperbolic terms. That the arithmetic of quadratic number fields arises is not surprising as 2-by-2 matrices have quadratic minimal polynomials. See e.g.\ \cite[Proposition 2]{St} for a general and completely explicit version of this formula. 

For the purpose of analytic number theory, the following variation of this procedure is very useful and leads to the Bruggeman-Kuznetsov formula \cite{Br, Ku}. Instead of integrating $K(x_1, x_2)$ over the diagonal $x_1 = x_2 \in X$, one can integrate it over $N(\Bbb{Z})\backslash N(\Bbb{R}) \times N(\Bbb{Z})\backslash N(\Bbb{R})$ (with $N = \{(\begin{smallmatrix} 1 & \ast \\ & 1\end{smallmatrix})\}$) against a   character   of  this group.  More precisely, we consider
$$\int_0^{\infty} \int_0^1\int_0^1 K\left(\xi_1 + i\frac{y}{|n|}, \xi_2 + i\frac{y}{|m|}\right) e(-\xi_1 n + \xi_2 m) \, {\rm d}\xi_1\, {\rm d}\xi_2  \frac{{\rm d}y}{y}$$
for $n, m \in \Bbb{Z}\setminus \{0\}$. Using the spectral expansion of $K$, the unipotent integrals feature the $n$-th and $m$-the Fourier coefficient of members of an orthonormal basis of $L^2(X)$. On the other hand, one can insert the definition \eqref{defK} of $K$ and compute the integral using the Bruhat decomposition of $\Gamma = {\rm SL}_2(\Bbb{Z})$, which in the present case can be stated very elementarily as
\begin{equation}\label{bruhat}
{\rm SL}_2(\Bbb{Z}) = N(\Bbb{Z}) \cup \bigcup_{c\in \Bbb{N}} \bigcup_{\substack{d \, (\text{mod } c)\\ (d, c) = 1}} N(\Bbb{Z}) \left(\begin{matrix} * & *\\ c & d\end{matrix}\right) N(\Bbb{Z}). 
\end{equation}
The details, due to Zagier, can be found in \cite[Section 1]{Jo}. This identifies the Kuznetsov formula as a relative trace formula in the sense of Jacquet. 
The final formula reads as follows:
\begin{equation}\label{kuz}
\begin{split}
&2\pi \sum_{\phi \in \mathcal{B}} \frac{\lambda_{\phi}(n) \lambda_{\phi}(m)}{L(1, \text{Ad}^2 \phi)} h(t_{\phi}) +   \int_{\Bbb{R}} \frac{\sigma_{it}(n) \sigma_{it}(m)}{|\zeta(1 + 2it)|^2} h(t)\,  {\rm d}t \\
& = \delta_{n, m} \int_{\Bbb{R}} h(t) \tanh(\pi t) t \,{\rm  d}t + \sum_c \frac{1}{c} S(n, m, c) \widehat{h}_{\pm}\left(\frac{nm}{c^2}\right)
\end{split}
\end{equation}
where $n, m \in \Bbb{Z} \setminus \{0\}$, $h$ is a sufficiently nice test function, $\mathcal{B}$ is an orthonomal basis of joint eigenfunction of $L^2_{\text{cusp}}(X)$, $\sigma_{it}(n) = \sum_{ad = |n|} (a/d)^{it}$ is the Hecke eigenvalue of the corresponding Eisenstein series, 
\begin{equation}\label{klooster}
S(n, m, c) = \sum_{\substack{d\, (\text{mod } c)\\ (d, c) = 1}} e\left( \frac{n d+ m\bar{d}}{c}\right)
\end{equation}
is the standard Kloosterman sum, and $\widehat{h}_{\pm}$ is a certain integral transform of $h$ given by a Bessel kernel depending on the sign $\pm = \text{sgn}(nm)$. The right hand side of \eqref{kuz} reflects directly the decomposition \eqref{bruhat}. 
The appearance of a Bessel function is not a coincidence, it is an archimedean analogue of the Kloosterman sum, as can be seen, for instance by the  formula \cite[8.432.7]{GR} 
$$K_0\left(\frac{\sqrt{nm}}{c}\right) = \frac{1}{2} \int_{\Bbb{R}} \exp\left( - \frac{1}{c} \left( n x + \frac{m}{x}\right)\right) \frac{{\rm d}x}{x}.$$
The adjoint $L$-function at 1 is  proportional to the square of the $L^2$-norm of a Hecke-normalized cusp form, and $|\zeta(1+2 it)|^2$ is the corresponding regularized version for Eisenstein series. 

There is a different way of proving the formula \eqref{kuz}, which is essentially Kuznetsov's original argument and analogous to the proof of the related Petersson formula for holomorphic cusp forms \cite[Proposition 14.5]{IK}.  The linear form $L^2(X) \rightarrow \Bbb{C}$, $\phi \mapsto \lambda_{\phi}(n)$ (say for $n > 0$) has a kernel  that is essentially given by a Poincar\'e series 
$$P_n(z) = \sum_{\gamma  \in N(\Bbb{Z}) \backslash \Gamma} f(n \cdot \Im \gamma z) e(n \cdot \Re \gamma z)$$
for a suitable function $f$. Indeed, by unfolding we have
$$\langle \phi,  P_n \rangle = \lambda_{\phi}(n) \sqrt{n}\int_0^{\infty} W_{it_{\phi}}(4\pi y) \overline{f(y)}  \frac{{\rm d}y}{y^2}.$$
Again by the Bruhat decomposition, one can compute the Fourier coefficients of $P_n$ in terms of Kloosterman sums. 
Computing now the inner product $\langle P_n, P_m\rangle$ by Parseval on the one hand and by unfolding on the other one derives \eqref{kuz}. It requires some analytic virtuosity with integral transforms and Bessel functions to play with the test function $f$ so that the left hand side of \eqref{kuz} features the given function $h$ and to compute the Bessel transform $\widehat{h}$. 

The shape of the Kuznetsov formula and the Selberg trace formula is similar in some respects, but there are also some important differences:
\begin{itemize}
\item The spectral side of the Kuznetsov formula does not contain the residual spectrum, in particular the discrete spectrum is cupsidal. 
\item The spectral side of the Kuznetsov formula is weighted by an $L$-value at the edge of the critical strip. For the group ${\rm GL}(2)$, such $L$-values can be removed and included with relatively little effort, see e.g. \cite[Proposition 26.15]{IK} for a prototype. For applications to $L$-functions involving period formulae it is often desirable to have an additional factor $1/L(1, \text{Ad}^2 \phi)$ in the cuspidal spectrum, but in other situations one may prefer a summation formula without an extra $L$-value. 
\item The Kloosterman sums in the Kuznetsov formula are easier to handle from the point of view of analytic number theory than the class numbers that appear in the Selberg trace formula. This is an important advantage of the Kuznetsov formula in practice (see the applications below). Theoretically, it is possible to switch directly between class numbers and sums of Kloosterman sums (cf.\ \cite[Theorem 1.3]{SY} or \cite{Al}). This is ultimately ``just'' Poisson summation, yet the passage is rather subtle. 
\end{itemize}
 
Unlike the Poisson summation formula, the Selberg trace formula and the Kuznetsov formula are asymmetric in the sense that the spectral side and the geometric side have a very different shape. It is an important feature that the Kuznetsov formula can nevertheless be inverted. For a sufficiently nice test function $H$ and integers $n, m \in \Bbb{Z} \setminus \{0\}$, the sum 
$$\sum_{c} \frac{1}{c} S(n, m, c) H\left( \frac{nm}{c^2}\right)$$
can be expressed in terms of Fourier coefficients of automorphic forms of the form
$$  \sum_{\phi \in \mathcal{B}} \frac{\lambda_{\phi}(n) \lambda_{\phi}(m)}{L(1, \text{Ad}^2 \phi)} \check{H}(t_{\phi})$$
for a suitable transform $ \check{H}$ of $H$ (essentially the inverse of $h \mapsto \hat{h}_{\pm})$, along with two similar sums involving the holomorphic spectrum and the Eisenstein series. 
It is the subtle inversion of the Sears-Titchmarsh transform $h \mapsto \widehat{h}_+$ (see \cite{ST}) that requires the \emph{entire} ${\rm GL}(2)$ spectrum including holomorphic cusp forms, see \cite[Section 16.4]{IK} for details. For the transform $h \mapsto \widehat{h}_-$, one needs the simpler Kontorovich-Lebedev inversion formula \cite[(16.46)]{IK}, where no holomorphic contribution shows up (which is consistent with the fact that holomorphic forms have no negative Fourier coefficients). 
  As we will see in the next section, the fact that \eqref{kuz} can also be read meaningfully from right to left is absolutely crucial for applications in analytic number theory and can hardly be underestimated. 
 
\section{Applications} 
 
 In this section we present a variety of applications of the Kuznetsov formula \eqref{kuz} and generalizations thereof to congruence subgroups, number fields and half-integral weight automorphic forms. They are roughly ordered in a thematic way, although there is no clear separation. Some applications use \eqref{kuz} from left to right, some from right to left  and some in both directions  usually with   a Voronoi step or a Cauchy-Schwarz step in between. Some of the problems are seemingly completely unrelated to the theory of automorphic forms and it is the sometimes accidental appearance of Kloosterman sums that opens the door. Needless to say that Kuznetsov formula is only one of many tools to solve many of the below mentioned problems, but it is a crucial one. Each of these problems deserves a survey article of its own, and I can only sketch a few ideas. I hope that the interested reader will consult the referenced literature for more details. 
  
 \subsection{Statistics of eigenvalues}\label{31}
 $\quad$
 
a) \emph{Vertical Sato-Tate laws.} Given a Hecke eigenform $\phi$, the statistical distribution of the eigenvalues $\lambda_{\phi}(p)$ as $p$ varies over primes has received a lot of attention. For holomorphic cusp forms  it is known to follow the Sato-Tate distribution \cite{BGHT}
$$\frac{\log P}{P} \sum_{\substack{p  \leq P\\ p \text{ prime}}} f(\lambda(p))\rightarrow \frac{1}{2\pi} \int_{-2}^2 f(x) \sqrt{4-x^2} \, {\rm d}x, \quad P \rightarrow \infty$$
($f$ a compactly supported continuous function),  but for Maa{\ss} forms this is wide open. One can modify the question, fix the prime $p$ and ask instead for the statistical distribution as one varies over the spectrum. This is known as a ``vertical Sato-Tate law''. It  was first addressed by Sarnak \cite{Sa} who showed using the Selberg trace formula that
$$\frac{12}{T^2} \sum_{\lambda_{\phi} \leq T} f(\lambda_{\phi}(p)) \rightarrow \frac{1}{2\pi} \int_{-2}^2 f(x) \sqrt{4-x^2} \frac{p+1}{p+2+\frac{1}{p} - x^2} \, {\rm d} x, \quad T \rightarrow \infty,$$
for a fixed prime $p$ and  $f$ as above. Serre \cite{Se} showed similar results for holomorphic forms using the Eichler trace formula. Note that if $p$ tends to infinity, this approaches the semicircle distribution. This question can also be addressed with the Kuznetsov formula (see e.g. \cite[Proposition 2]{BBR}), in which case the spectral sum is naturally weighted by an $L$-value:
$$\frac{12}{T^2} \sum_{\lambda_{\phi} \leq T} f(\lambda_{\phi}(p)) \frac{\zeta(2)}{L(1, \text{Ad}^2 \phi)} \rightarrow \frac{1}{2\pi} \int_{-2}^2 f(x) \sqrt{4-x^2}  \, {\rm d}x, \quad T \rightarrow \infty.$$
Interestingly, this slightly different counting procedure produces the semicircle distribution ``on the nose''. For the proof, we assume by a standard approximation argument  that $f(x) = x^k$ for some $k \in \Bbb{N}$, apply the Hecke relations and estimate  the Kloosterman term in the Kuznetsov formula for instance by the Weil bound \eqref{we} below. \\

b) \emph{Density results for Maa{\ss} forms violating the Ramanujan conjecture.} The previous two asymptotics  show in particular that the Ramanujan conjecture at $p$ is satisfied ``with probability one''.  By a small variation of the argument, this can be made quantitative as follows (\cite[Proposition 1]{BBR}): for a prime $p \leq T$, $\alpha > 2$ and $\varepsilon > 0$ we have
\begin{equation}\label{dens}
 \# \{ \lambda_{\phi} \leq T^2 : |\lambda_{\phi}(p)| \geq \alpha\} \ll_{\varepsilon} T^{2-\frac{8 \log (\alpha/2)}{\log p} + \varepsilon}.
\end{equation}
Upon choosing $\alpha = p^{1/4 + \varepsilon}$, this recovers in particular the Selberg-type bound $\lambda_{\phi}(p) \ll p^{1/4 + \varepsilon}$ (since the right hand side of \eqref{dens} is less than 1, so the set on the left hand side must be empty). It is at this point where the Kuznetsov formula is advantageous compared to the Selberg trace formula \cite[Theorem 1]{Sa}. The latter would replace the constant 8 in the exponent only by 4,\footnote{which follows from Sarnak's argument \cite{Sa} with small modifications} and the corresponding version of \eqref{dens} would only interpolate between the trivial representation (with $\lambda_{\text{triv}}(p) = p^{1/2} + p^{-1/2}$) and the tempered spectrum, that is, the bound is trivial on the tempered spectrum $\alpha = 2$ and excludes the possibility of $\lambda_{\phi}(p)|\geq p^{1/2 + \varepsilon}$. Had we used trivial bounds on the Kloosterman sum, the Selberg trace formula and Kuznetsov formula would produce the same result, but the option of getting extra cancellation from the Weil bound  \eqref{we} makes the Kuznetsov formula advantageous in this situation.   

This type of argument works also at the infinite place. Using the test function $h(t) = (X^{it} + X^{-it})^2(t^2 + 1)^{-2}$ that blows up at the exceptional spectrum, one obtains Selberg's $3/16$ bound as a lower bound for non-trivial eigenvalues of the Laplacian on $\Gamma \backslash \Bbb{H}$, see \cite[(16.58)]{IK} as well as density results for exceptional eigenvalues of large level. \\

c) \emph{Large sieve estimates.} A different kind of statistical result is provided by large sieve estimates that show in some sense orthogonality of Hecke eigenvalues against arbitrary families of sequences. More precisely, if   $a_n$ is any sequence of complex numbers, then \cite[Theorem 2]{DI}
$$\sum_{t_{\phi} \leq T} \Big|\sum_{n \leq N}a_n \lambda_{\phi}(n) \Big|^2 \ll_{\varepsilon} (T^2 + N)^{1+\varepsilon} \sum_{n \leq N} |a_n|^2$$
 for all $\varepsilon > 0$. This is essentially best possible, since the first term on the right hand side is achieved when all $n$-sums have square-root cancellation and the second term is achieved if $a_n$ is chosen to be $\overline{\lambda_{\phi_0}(n)}$ for some fixed $\phi_0$. 
 Note that a simple application of Cauchy's inequality would replace $T^2 + N$ with $T^2N$. Such general estimates are very useful in analytic number theory. One of the most common applications is a bound for moments of $L$-functions of a family $\mathcal{F}$ that recovers in many cases the Lind\"of hypothesis on average over the family if the conductor of the $L$-functions is at most $|\mathcal{F}|^4$ (this follows from an approximate functional equation). The proof of the large sieve inequality is a double application of the Kuznetsov formula. Opening the square, one applies the Kuznetsov formula to the $t_{\phi}$-sum. After an application of Cauchy's inequality, one can re-arrange the Kloosterman term to make it amenable to a second application of the Kuznetsov formula in the other direction. This is not involutory because of Cauchy's inequality in between. 
 
 \subsection{Arithmetic applications: Kloosterman sums, primes and arithmetic functions}\label{33}
 
$\quad$
 
a) \emph{Linnik's problem: cancellation in sums of Kloosterman sums.} Kloosterman sums as in \eqref{klooster} are ubiquitous in number theory. They first came up  in Kloosterman's investigation \cite{Klo} of quaternary quadratic forms by the circle method. Weil's bound  \cite{We}
\begin{equation}\label{we}
|S(m, n, c)| \leq \tau(c) (m, n, c)^{1/2} c^{1/2}
\end{equation}
(where $\tau$ denotes the divisor function) is best possible for individual sums, but quite often in number theory \emph{sums of Kloosterman sums} 
$$\mathcal{S}_{m, n}(X; W) := \sum_{c} S(n, m, c) W(c/X)$$
over the modulus $c$ come up for some fixed smooth function $W$ and some large parameter $X$. Weil's bound yields immediately $\mathcal{S}_{m, n}(X; W) \ll_{W, m, n} X^{3/2} \log X$, but the Kuznetsov formula (along with the fact that ${\rm SL}_2(\Bbb{Z})$ has no exceptional eigenvalues, see \cite[p.\ 261]{DI} for a quick proof) provides substantial cancellation in this sum: 
\begin{equation}\label{Linnik}
  \mathcal{S}_{m, n}(X; W) \ll_{W, m, n, \varepsilon} X^{1+\varepsilon}
  \end{equation} for all $\varepsilon > 0$. This is very remarkable because it goes far beyond what algebraic geometry can achieve and it is one of the very few examples 
  where exponential sums can be averaged non-trivially over the modulus. 
  The bound \eqref{Linnik} and generalizations thereof, essentially due to Kuznetsov \cite{Ku} and developed in \cite{DI}, were historically the first applications of \eqref{kuz} and the cornerstone for all applications in this subsection. \\

b) \emph{Shifted convolution sums.} A typical question in number theory concerns the correlation of additive shifts of multiplicative function. A standard example is an asymptotic evaluation of
$$\sum_{n \leq x} \tau(n) \tau(n+1)  = \sum_{a \leq x} \sum_{\substack{m \leq x\\ m \equiv 1 \, (\text{mod } a)}}\tau(m)$$
(as can be seen by opening the first divisor function and writing $n = ab$). The inner sum can be transformed by the Voronoi summation formula into a sum over Kloosterman sums $S(m, 1, a)$. It is then the Kuznetsov formula that can evaluate the $a$-sum over Kloosterman sums with respect to the modulus very precisely and yields the strongest error terms.\footnote{It is also possible to detect the condition additive shift by some variant of the circle method. After some elementary Fourier analysis this also leads to Kloosterman sums. Although relatively similar, in many cases the above approach has slight advantages.}  The (somewhat involved) details can for instance be found in \cite{Me}. If the sum over $x$ is smooth, one obtains square root cancellation in the error term. Such estimates and generalizations thereof play an important role in the theory of the Riemann zeta-function (see below). 
\\

c) \emph{Primes in long arithmetic progressions.} The distribution of prime numbers 
$$\pi(x; q, a) := \#\{p \leq x\mid p \equiv a \, (\text{mod } q), p \text{ prime}\}$$
in arithmetic progressions with $(a, q) = 1$ is a very delicate topic due to possible Siegel zeros of Dirichlet $L$-functions. For many applications one can get away with a few arithmetic progressions in which primes are potentially ill-distributed if the majority of progressions is well-behaved. The standard result of this type is the Bombieri-Vinogradov theorem (\cite[Theorem 17.1]{IK}) which essentially states that $\pi(x; q, a) \sim  \phi(q)^{-1} \pi(x; 1, 0)$ for most pairs $a, q$ with $(a, q) = 1$ as long as $q \leq x^c$ for some $c < 1/2$. Going beyond $c > 1/2$ can be crucial in applications. In a sequence of papers culminating in \cite{BFI}, Bombieri, Friedlander, Fouvry and Iwaniec managed to achieve this at least in certain situations, a cruical ingredient being non-trivial bounds for sums over (incomplete) Kloosterman sums as they follow from the Kuznetsov formula.  Although these extended versions of the Bombieri-Vinogradov theorem  have some restrictive conditions (e.g.\ $a$ has to be fixed), it is remarkable that they go beyond the ranges that would be covered by the Riemann hypothesis. This can hardly be underestimated and is the basis of the following application.\\

d) \emph{The first case of Fermat's last theorem.} If $p$ is a prime, then the first case of Fermat's last theorem gives insolubility of the equation $x^p + y^p + z^p = 0$ in integers $x, y, z$ such that $p \nmid xyz$. Adleman and Heath-Brown \cite{AHB} showed that this follows for infinitely many primes (in a strong quantitative sense) if a certain estimate on primes in arithmetic progression holds; this was proved by Fourvy \cite{Fo} based (among other things) on the work of Bombieri-Friedlander-Iwaniec \cite{BFI}. \\

 e) \emph{Proportion of zeta zeros on the critical line.} Selberg proved in 1942 that a positive proportion of zeros lie on the critical line $\Re s = 1/2$. About 30 years later, Levinson showed   by a different method that in fact at least 1/3 of all zeros lie on the critical line (see   \cite[Section 10.1]{Iv} for a brief sketch of the method). It ultimately depends on computing the mean square of a linear combination of $\zeta(s)$ and its derivatives multiplied by a carefully chosen Dirichlet polynomial whose length is responsible for the numerical value of the proportion. Based on mean values for (incomplete) Kloosterman sums following from \eqref{kuz}, Conrey \cite{Co}  improved the proportion of zeros on the critical line to slightly over 40\%.\footnote{At the time of writing, the current record is slightly over 41\%.} \\

f) \emph{The largest prime factor of $n^2+1$.} It is an old conjecture that there are infinitely many primes of the form $n^2+1$, but a proof is unfortunately far out of reach by current methods. As an approximation, one can consider the largest prime factor of $n^2+1$. It is trivial that this can be at least as big as $n$ infinitely often\footnote{For any prime $p \equiv 1$ (mod 4) choose $1 \leq n < p$ with $n^2 \equiv -1$ (mod $p$).}. Chebychev indicated a method to show that it is bigger than $Cn$ infinitely often for every given constant $C > 1$. This was refined by Hooley. The key point is that the distribution of $n^2+1$ in residue classes leads to Kloosterman sums (cf.\ \cite[Lemma 2]{DI2}): one has 
$$\sum_{n^2 + 1 \equiv 0 \, (\text{mod } m)} e\left(\frac{nh}{m}\right) \approx \sum_{\substack{r^2 + s^2 = m\\ r, s > 0, (r, s) = 1}} e\left( \frac{h\bar{r}}{s}\right)$$
which after Poisson summation in $r$ gives Kloosterman sums. Using bounds on sums of Kloosterman sums, Deshouillers and Iwaniec showed \cite{DI2} that the largest prime factor of $n^2+1$ is infinitely often at least $n^{6/5}$. At the time of writing, the current record is $n^{1.279}$ \cite{Mer}.
\\ 
 
g) \emph{Equidistribution of roots of quadratic congruences.} The previous example is ultimately concerned with the equidistribution of roots of quadratic congruences, a topics that goes back (in principle to Gau{\ss} and) Hooley, Hejhal, Bykovski\u{i}, Zavorotny and others with the strongest results based on the spectral theory of automorphic forms. Duke, Friedlander and Iwaniec \cite{DFI} went one step further and coupled this with a sieve to show that the equidistribution property remains true among \emph{prime} moduli: let $f$ be an integral quadratic polynomial without real roots and  $h \not= 0$. Then
$$\sum_{\substack{p \leq x\\ p \text{ prime}}} \sum_{\substack{\nu \, (\text{mod } p)\\ p \mid f(\nu)}} e\left(\frac{h\nu}{p}\right) = o\Big(\sum_{p \leq x} 1\Big)$$
as $x \rightarrow \infty$. A beautiful set of lecture notes proving this result \emph{ab ovo} is \cite{Ko}. The case of polynomials with real roots is treated in \cite{To}. \\
 
h) \emph{Prime geodesic theorem.} The prime geodesics on the upper half plane play a similar role for the Selberg zeta functions as the primes for the Riemann zeta function. In particular they satisfy a ``prime number theorem'', which follows from a standard application of the Selberg trace formula: if $p$ runs through the lengths of prime geodesics on $\Gamma \backslash \Bbb{H}$, then
$$\sum_{ p \leq X } \log p = X + O(X^{3/4}).$$
Henryk Iwaniec observed that an \emph{additional} application of the Kuznetsov formula can improve the error term. These ideas have been refined, and the current record \cite{SY} is an error term of $O(X^{25/36+\varepsilon})$ with $25/36 = 0.69\ldots$ which uses the Kuznetsov formula again, but rather indirectly through the work of Conrey-Iwaniec \cite{CI}.

 \subsection{Applications to $L$-functions I}\label{32}
 
$L$-function occur naturally in families. While their individual behaviour is often hard to understand, their statistical properties in families are more amenable to the tools of analytic number theory. If the family is given by spectral properties, spectral summation formulae such as the Kuznetsov formula belong to the key tools to study analytic properties of $L$-functions in various average senses. Therefore all of the following applications of the Kuznetsov formula are concerned with moments of $L$-functions in one way or another. Usually the $L$-functions are encoded into the Kuznetsov formula by means of an approximate functional equation, but in principle it is also possible to use a suitable integral representation.    This area is a huge industry, and the following discussion can highlight, inevitably, only a very small selection of results. \\

 a) \emph{Symmetry types of $L$-functions.} In \cite{KS}, Katz and Sarnak introduced the notion of symmetry type for a family of $L$-functions. A common approach to this determination  is to analyze the density of low-lying zeros in the specified family which should mimic the appropriate random matrix model. With current technology this is only possible for test functions whose Fourier transforms have restricted support. The most frequently used quantity is the so-called one-level density which essentially counts the number of zeros within $1/\log C$ of the origin (where $C$ is the conductor of the $L$-function). The key tools here are a spectral summation formula (like the Kuznetsov formula or the Petersson formula for holomorphic forms) and Weil's explicit formula (see e.g.\ \cite[Theorem 5.11]{IK}) to convert zeros into primes, and the results typically assume the Riemann hypothesis.  One of the most influential papers in this direction  is \cite{ILS}, covering many families of $L$-functions. The one-level density is not unrelated to the analytic rank of the $L$-function (i.e.\ its order at $s=1$). Applications to ranks of elliptic curves are given for instance in \cite{Yo}.   \\
 
 b) \emph{Spectral expansion of the fourth moment of the Riemann zeta function.} Motohashi observed that for  a Schwartz class function $w$, the fourth moment
 $$\int_{-\infty}^{\infty} |\zeta(1/2 + it)|^4 w(t)\, {\rm d}t$$
 of the Riemann zeta function has a beautiful expansion as a spectral sum of third powers of central $L$-values of automorphic forms for ${\rm SL}(2, \Bbb{Z})$, see \cite{Mo} for a detailed and self-contained proof and extensive discussion. This formula is quite remarkable from an aesthetic point of view, relating two rather different families of $L$-functions in an exact identity, but it is also very useful theoretically and gives the strongest bounds for the fourth moment in short intervals (by choosing $w$ appropriately and analyzing the dual side of the formula). The Kuznetsov formula enters after opening the fourth moment: using an approximate functional equation, we obtain an expression very roughly of the form
 $$\int_{-\infty}^{\infty} \sum_{n, m \ll (1+ |t|)^{1/2}} \frac{\tau(n)\tau(m)}{n^{1/2 + it}m^{1/2 - it}} w(t) \, {\rm d}t. $$
 Here it is convenient to sum the diagonals and write $m = n+h$, so that shifted convolution sums of the type  $\sum_n \tau(n) \tau(n+h)$ (with some weights) emerge.
 \\
 
c) \emph{Beyond endoscopy.} Given a reductive algebraic group $G$, the cuspidal automorphic representations $\pi$ of $G(\Bbb{A})$  are expected to be functorial transfers from other groups if $L(s, \pi, \rho)$ has a pole for some finite dimensional representation $\rho: {\,} ^LG \rightarrow {\rm GL}_n(\Bbb{C})$. Langlands \cite{La} put forward the idea to study this by the trace formula. The idea is to consider a trace formula whose spectral side  is weighted by $\text{ord}_{s=1} L(s, \pi, \rho)$ and to try to match the geometric side with the geometric side of the corresponding other group.  Sarnak suggested  to study this with the  Kuznetsov formula instead.  In \cite{Ve}, Venkatesh carried out this programme with the Kuznetsov formula for the case $G = {\rm GL}(2)$ and $\rho$ the symmetric square representation, which indeed classifies precisely the dihedral representations. It would be very interesting to study higher rank groups or higher dimensional representations $\rho$, and we refer the reader also to \cite{Al} where Langlands' ideas of  ``Beyond endoscopy'' are investigated with the Selberg trace formula. \\ 

d) \emph{Non-vanishing of $L$-functions.} Evaluating asymptotically a moment of a family of $L$-functions can prove that one, infinitely many, or even a large proportion of members in this family are non-zero. Whether or not an $L$-function vanishes at a certain point in the critical strip (often the central point $s=1/2$) is an important question, sometimes with unexpected applications. 

If $\Gamma \subseteq {\rm SL}_2(\Bbb{Z})$ is a finite-volume, but non-cocompact subgroup, all but finitely many eigenvalues are embedded in the continuous spectrum, and it is in general unclear  how to prove a Weyl law for the discrete spectrum alone. In fact, it may well be that there is no discrete spectrum at all except for the constant function.   Philipps and Sarnak  \cite{PS} considered deformations   given by curves in the Teichm\"uller space of $\Gamma$, and for congruence subgroups $\Gamma$ they gave a  criterion in terms of non-vanishing of $L$-functions which tells us (conditionally under the assumption of the simplicity of the spectrum of $L^2({\rm SL}_2(\Bbb{Z}) \backslash \Bbb{H})$) that such deformations destroy cusp forms,  i.e.\ the eigenvalue becomes a pole of the scattering matrix\footnote{This can be phrased as an analogue of Fermi's golden rule.}. Using the Kuznetsov formula, Luo \cite{Lu},  evaluated  a mollified second moment of the respective Rankin-Selberg $L$-functions to show that indeed a positive proportion of these $L$-values doesn't vanish.  
 
Also higher rank examples can be treated with Kuznetsov formula \eqref{kuz}: in \cite{BLM} it is shown that   for any self-dual cuspidal automorphic representation $\Pi$ on ${\rm GL}(4, \Bbb{Q}) \backslash {\rm GL}(4, \Bbb{A})$ which is unramified at all finite places, there exist infinitely cusp forms $\pi$ for ${\rm SL}(2, \Bbb{Z})$ such that $L(1/2, \Pi \times \pi) \not= 0$. This is the consequence of a certain reciprocity formula relating different families of $L$-functions in an exact identity.

\subsection{Applications to $L$-functions II: subconvexity and equidistribution}\label{34}
  
The Generalized Riemann Hypothesis implies best possible bounds for $L$-functions on the critical line (cf.\ \cite[Corollary 5.20]{IK}). In absence of a proof of GRH, it is a very interesting and very hard problem to obtain bounds for $L$-functions that are better than the generic convexity bound that is implied by the functional equation. Such bounds can often imply that ``something or other is equidistributed'' \cite[p. 373]{Fr}, and we give a number of examples where the corresponding subconvexity bound is a consequence of the Kuznetsov formula. From a technical and conceptual point of view, these are often the hardest applications of the Kuznetsov formula, since the relative trace formula is coupled with various other transformation formulae (Voronoi, Poisson, or even another application of the Kuznetsov formula), so that  extremely precise knowledge on the exact shape of the various terms in the formula is required. \\

a)  \emph{Sums of three squares.} An asymptotic formula for the number of solutions to $x^2 + y^2 + z^2 = n$ should  in principle follow from Siegel's mass formula: the ternary theta series is decomposed into an Eisenstein part that gives the main term as a product of local densities, and a cuspidal part, whose Fourier coefficients should be substantially dominated by the main term and therefore constitute an error term. There are many subtleties for positive ternary quadratic forms (see e.g.\ \cite{Bl} for a survey), but a major problem is our limited knowledge towards the Ramanujan conjecture for half-integral weight cusp forms. The key breakthrough here is due to Iwaniec and Duke \cite{Iw, Du}, based on analyzing the Sali\'e-type sums in the half-integral weight Kuznetsov formula. There is an alternative route based on Waldspurger's theorem that translates the problem into a problem of obtaining subconvexity for twisted automorphic $L$-functions. This can also be achieved by Kuznetsov formula (applied in a very different way) and gives the currently strongest results \cite{BH}. Once an asymptotic formula with a power saving error is available, one can ask the finer question of equidistribution of the solutions in shrinking  regions on the sphere of radius $\sqrt{n}$ which was first obtained in \cite{DSP}. A striking result is an analysis of the variance in small annuli \cite{HR}, which also uses the Kuznetsov formula. \\

b) \emph{Equidistribution of Heegner points.}  A similar problem asks for equidstribution of Heegner points and generalizations thereof which can again be rephrased as a subconvexity problem of twisted $L$-functions. Here the quadratic form $x^2 + y^2 + z^2$ is replaced with the discriminant form $q(x, y, z) = y^2 - 4xz$. For negative $d$, one associates to each ${\rm SL}_2(\Bbb{Z})$-orbit of a primitive solution to $q(x, y, z) = d$ the ${\rm SL}_2(\Bbb{Z})$-orbit of the Heegner point $(- y + i\sqrt{|d|})/(2x) \in \Bbb{H}$. After a conditional proof of Linnik, Duke showed in his seminal paper   \cite{Du} that these points equidistribute on the modular surface. This has initiated  a lot of work, and we refer the reader for instance to \cite{Mi, HM, Coh, Zh, LMY, Yo2} for generalizations and extensions. See \cite{MV} for an excellent survey.

Duke's theorem has been generalized to higher rank and in particular to ${\rm GL}(3)$ in \cite{ELMV} where the authors obtain analogous equidistribution result for equivalence classes of maximal compact flats. As in Duke's original  theorem, this can be rephrased in algebraic terms and gives statements about the distribution of orders in totally real cubic fields and of integer points on certain homogeneous varieties. For instance,  if $P_i$ runs through a sequence of integral, cubic, monic, irreducible polynomials with three real roots and increasing discriminant, then (a projection of) the set of integral $3 \times 3$ matrices $M$ with $P_i(M)=0$  becomes equidistributed with respect to a natural ${\rm PGL}_3(\Bbb{R})$-invariant measure. One (of many) inputs is a subconvexity bound for automorphic $L$-functions with nebentypus established in \cite{BHM} using the ${\rm GL}(2)$ Kuznetsov formula. \\

c) \emph{Mass distribution problems.}   A different class of problems is concerned with the equidistribution of the mass of a sequence of automorphic forms whose spectral parameter (or potentially also its level) is tending to infinity. The random wave conjecture  of Berry \cite{Be} states that eigenfunctions of a classically ergodic system behaves like Gaussian random variables in the large eigenvalue limit. Although some exceptions are known \cite{RS}, the mass of (a sequence of) automorphic forms $\phi_j$, $j = 1, 2, \ldots \rightarrow \infty$,  is generically expected to become equidistributed. This can be measured in at least three interesting ways:
\begin{itemize}
\item in the weak-$\ast$ limit: $|\phi_j|^2 {\rm d}\mu \rightarrow c \, {\rm d}\mu$ for a  constant $c > 0$. This is related to the QUE conjecture of Rudnick and Sarnak \cite{RS}.  For Hecke eigenfunctions on ${\rm SL}_2(\Bbb{Z})\backslash \Bbb{H}$ this is related to subconvexity of certain degree 8 $L$-functions (degree 4 $L$-functions for Eisenstein series). For Eisenstein series and cuspidal CM-forms (dihedral forms), much progress has been made for instance in \cite{PeS, Sa3, LS} long before the breakthrough results of Lindenstrauss and Soundararajan. More refined quantitative recent results can be found for instance in \cite{Hu, BK2, Yo3}. All of them are based on certain mean value bounds of $L$-functions by the Kuznetsov formula. 
\item  by considering $\| \phi_j \|_p/ \| \phi_j \|_2$ for some $2 < p \leq \infty$. Bounds for this quantity show that large peaks cannot occur, or at least can occur only rarely. Particularly interesting is the case $p=4$, because by Parseval and the formula of Watson-Ichino \cite{Ic} this can be expressed as a mean value of central $L$-values, see for instance  \cite{Bl2, BKY, BK, Hu}.  
\item by considering restriction norm problems to interesting submanifolds $Y \subseteq X$. In some interesting cases the restriction $\| \phi|_Y \|^2_{2, Y}/\| \phi \|^2_{2, X}$ can be expressed in terms of central values of $L$-functions. Typical situations are Gross-Prasad pairs (see \cite{II}) of subgroups, but also the period integral for ${\rm GL}_n \times {\rm GL}_{n-1}$ Rankin-Selberg convolutions belongs into this class. By Parseval this can often be translated into a mean value of $L$-functions that in favourable cases can be attacked with the Kuznetsov formula, see e.g.\ \cite{BKY, LY, LLY} and \cite{Sa-letter} for a general overview. 
\end{itemize}
 
Finally it is worth mentioning that subconvexity of $L$-functions has become an important topic in the theory of $L$-functions in its own sake. On the one hand it is an excellent test case of the strength of available methods, on the other hand the search for subconvexity has also led to a better understanding of available tools and in particular the Kuznetsov  formula. Even in absence of arithmetic applications it broadens our understanding of $L$-functions. This is huge area, and we can only mention here a very small set of selected examples focusing on higher rank $L$-functions, for instance \cite{Bl10, JM, Li, Su}.

\section{Other groups} 

After the great success of the Kuznetsov formula for the group ${\rm GL}(2)$, it is natural to investigate what it has to offer in other situations. The generalization to other rank 1 groups does not pose essential difficulties, and the shape of the formula is relatively similar. In fact, the formalism presented in Section \ref{sec2} works for any reductive algebraic group. General versions for ${\rm GL}(n)$ can be found in  \cite[Section 11]{Go} and \cite[Theorem D]{Frb}. Applications in analytic number theory require very precise and explicit information on both sides of the formula, and this is not at all easy to obtain for higher rank groups. In particular, some of the key open questions involve 
\begin{itemize}
\item a good understanding of the various types of Kloosterman sums as certain kinds of exponential sums over finite fields and rings;
\item a good understanding of the growth and oscillatory properties of the relevant integral transforms, and also the problem to what extent they can be inverted. 
\end{itemize}
Nevertheless, the Kuznetsov formula is a versatile tool for statistical investigations of automorphic forms and $L$-functions potentially for arbitrary groups and perhaps analytically easier to handle than the Selberg trace formula. A thorough study of its scope and limitations in higher rank has just begun, but is a rapidly growing body of work. The following discussion highlights a few specific examples.

\subsection{The group ${\rm GL}(3)$}
 
The Weyl group of ${\rm GL}(3)$ has 6 elements, so the right hand side of the Kuznetsov formula consists (a priori) of 6 terms. Two of them vanish, one is the identity, and we are left with two hyper-Kloosterman terms and the long Weyl element Kloosterman term. The latter is usually dominant in applications.  The first explicit computations go back to the late 1980's in \cite{BFG} where all relevant Kloosterman sums were completely determined, and their algebro-geometric properties were studied in \cite{Frb, Ste}. To get a feeling for the complexity, an explicit version of the long Weyl element Kloosterman sum for ${\rm GL}(3)$ looks as follows: let $n , m  \in( \Bbb{Z} \setminus\{0\})^2$ and $\psi_{n}$, $\psi_{m}$ be the corresponding characters on unipotent upper-triangular $3 \times 3$ matrices, explicitly $$\psi_{n}\left(\begin{smallmatrix} 1 & x_2 & \ast \\ & 1 & x_1 \\ & & 1\end{smallmatrix}\right) = e(n_1x_1 + n_2 x_2).$$
Let $w$ the long Weyl element, $c =\text{diag} (1/c_2, c_2/c_1, c_1)$ for some $c_1, c_2\in\Bbb{N}$ 
 and $U$ the integral unipotent upper-triangular $3\times 3$ matrices.    Then by the explicit Bruhat decomposition the long Weyl element Kloosterman sum corresponding to the ``modulus'' $c$ and the characters $\psi_n$, $\psi_m$ is given by 
 $$S(\psi_n, \psi_m, c) = \sum_{\gamma \in b c w b' \in U\backslash{\rm SL}_3(\Bbb{Z})/ U} \psi_n(b) \psi_m(b').$$
Unfolding the definition, this equals
\begin{equation}\label{weyl}
 \sum_{\substack{B_1, C_1 \, ({\rm mod }\, c_1)\\B_2, C_2 \, ({\rm mod }\,  c_2)\\ c_1C_2 + B_1B_2 + c_2C_1 \equiv 0 \, ({\rm mod }\, c_1c_2)\\ (B_j, C_j, c_j) = 1}} \hspace{-0.7cm} e\left(\frac{m_2B_1 + n_1(Y_1 c_2 - Z_1 B_2)}{c_1} - \frac{m_1B_2 + n_2(Y_2 c_1 - Z_2B_1)}{c_2}\right)
 \end{equation}
where $Y_j, Z_j$ are chosen such that $Y_jB_j + Z_jC_j \equiv 1$ (mod $c_j$) for $j=1, 2$. If $(c_1, c_2) = 1$, this factorizes into a product of two ordinary Kloosterman sums to modulus $c_1$ and $c_2$.  The archimedean Whittaker transform is even more complicated and only recently has been made explicit by Buttcane \cite{Bu1}. 

With these developments in place, many of the applications in Section \ref{31} can be performed in a similar way (although sometimes quantitatively weaker): vertical Sato-Tate laws and density results for exceptional eigenvalues can be found along these lines in  \cite{BBR, Bl5,  BBM,  BZ}. A large sieve has been established in \cite{BB1}. 

The invertibility of the Kuznetsov formula was the key to all applications in Section \ref{33}. First important steps in this direction were established in   \cite{Ye, Bu2}, the most complete solution is contained in \cite{Bu3}, where a general test function can be put on the long Weyl element term. As one may expect, one of the difficulties lies in the fact that the spectral side contains not only the spherical spectrum, but all weights simultaneously. In the case of a non-abelian maximal compact subgroup $K = {\rm SO}(3)$, this is a non-trivial obstacle. However,  it remains an open question where (and if) the long Weyl element Kloosterman sums \eqref{weyl} come up ``in nature'', i.e.\ in number theoretic problems.  Applications of arithmetic significance are still a desideratum.

Some work has been done on $L$-functions corresponding to the applications presented in Sections \ref{32} and \ref{34}. Symmetry types of $L$-function using the one-level density have been investigated in \cite{GK}. The underlying analysis of the Kuznetsov formula was recently generalized to ${\rm GL}(4)$ in \cite{GSW}.

The first  subconvexity result for genuine ${\rm GL}(3)$ $L$-functions (by any method) has been obtained in \cite{BB2} using the Kuznetsov formula: let $\pi   \subseteq L^2({\rm SL}_3(\Bbb{Z})\backslash {\rm SL}_3(\Bbb{R}))$ be an irreducible, cuspidal, spherical representation with  Langlands  parameter $\mu$   in ``generic position'', i.e.\ satisfying  
$$c \leq \frac{|\mu_{ j}|}{\| \mu\|} \leq C \quad (1 \leq j \leq 3) \quad \text{and} \quad c \leq \frac{|\mu_{ i} - \mu_{ j}|}{\| \mu\|} \leq C \quad (1 \leq i< j \leq 3)$$
for two constants $C > c > 0$. This covers $99\%$ of all Maa{\ss} forms (choosing $c, C$ appropriately), but misses, for instance, self-dual forms. 
Then 
 $$L(1/2, \pi ) \ll  \| \mu \|^{\frac{3}{4} - \frac{1}{120000}},$$
where the convexity bound is $L(1/2, \pi ) \ll  \| \mu \|^{3/4 + \varepsilon}$.

\subsection{The group ${\rm GSp}(4)$}

Another interesting group  with real rank 2 is the group ${\rm Sp}(4)$. As the Siegel upper half plane admits a complex structure, we have the notion of holomorphic cusp forms, for which an analogue of the Petersson formula has been worked out in detail by Kitaoka \cite{Ki}. This formula has recently been used in analytic number theory. For Siegel modular forms, there is a big difference between Fourier coefficients and Hecke eigenvalues. The unipotent periods produce Fourier coefficients in the relative trace formula, and it is a remarkable fact that by B\"ocherer's conjecture \cite{Bo, FM} they are related to central $L$-values. 

Using the Petersson-Kitaoka relative trace formula, local spectral equidistribution problems in the sense of Section \ref{31} have been considered in \cite{KST}. The papers \cite{Bl4, Wa} compute mean values of spinor $L$-functions with applications to non-vanishing. This uses the full force of the formula. The analysis of the proof   manipulates the Kloosterman terms in a rather sutble way, as they also contribute to the main term in the asymptotic formula.

\subsection{Groups of unbounded rank} The asymptotic distribution of Satake parameters has been obtained in a rather general context.  
Important works in this direction using the   currently most advanced analysis of the Selberg trace formula for the group ${\rm GL}(n)$ and even arbitrary classical groups are \cite{MT, FiMa}. Using a version of the Kuznetsov formula, strong density results for Satake parameters are obtained in \cite{Bl5} for the group
$\Gamma_0(q) \subseteq {\rm SL}_n(\Bbb{Z})$ of matrices whose lowest row is congruent to $(0, \ldots, 0, \ast)$ modulo $q$. This is based, among other things, on a detailed analysis of certain ${\rm GL}(n)$ Kloosterman sums.  For a discussion of these results in particular in connection with Sarnak's density conjecture we refer the reader to the introduction of \cite{Bl5}. A Petersson type formula is used in \cite{KL} to obtain asymptotic distribution results for Satake parameters on ${\rm GSp}(2n)$.

\end{document}